\documentclass[11pt,a4paper]{article}
\usepackage{amsfonts,amsmath,amssymb}
\usepackage{mathrsfs}
\usepackage[latin1]{inputenc}
\usepackage{mathrsfs}
\usepackage[latin1]{inputenc}
\usepackage{amsmath,amsfonts,amsthm}

\usepackage{mathrsfs}

\usepackage{times}
\usepackage[latin1]{inputenc}

\usepackage{mathrsfs}

\newtheorem{theo}{Theorem}[section]
\newenvironment{Theo}{\begin{theo}\slshape}{\end{theo}}
\newtheorem{lemm}{Lemma}[section]
\newenvironment{Lemm}{\begin{lemm}\slshape}{\end{lemm}}

\newtheorem{Corol}{Corollary}[section]
\newtheorem{rema}{Remark}[section]

\newtheorem{exam}{Example}[section]

\newtheorem{defi}{Definition}[section]

\newcommand{\tild}{\widetilde}
\newcommand{\p}{\partial}

\def\qed{\hfill$\square$\par \bigskip}
\newenvironment{Demo}[1]{{\bf Proof#1.~}}{\qed}
\newcommand{\N}{\mathbb{N}}
\newcommand{\R}{\mathbb{R}}

\newcommand{\dive}{\mathrm{div\,}}
\newcommand{\curl}{\mathrm{curl\,}}
\newcommand{\Ldiv}{H(\dive;\,\Omega)}

\newcommand{\uu}{\bold u}

\newcommand{\vvv}{\bold v}
\newcommand{\www}{\bold w}

\newcommand{\va}{\varphi}
\newcommand{\ep}{\varepsilon}

\newcommand{\ddd}{\textrm{div}\thinspace}
\newcommand{\DDelta}{\Delta_{\mu,\lambda}}

\newcommand{\norm}[1]{\|#1\|}
\newcommand{\abs}[1]{\left\vert#1\right\vert}
\newcommand{\set}[1]{\left\{#1\right\}}
\newcommand{\para}[1]{\left(#1\right)}
\newcommand{\cro}[1]{\left[#1\right]}

\newcommand{\seq}[1]{\left<#1\right>}

\newcommand{\To}{\longrightarrow}
\newcommand{\NN}{\mathcal{N}_{\tau,\varphi}(\vvv)}
\newcommand{\NNN}{\mathcal{N}_{\tau,\varphi}(\widetilde{\vvv})}

\begin{document}
\title{\bf \Large Carleman Estimate and Inverse Source
Problem for Biot's Equations Describing Wave Propagation
in Porous Media}
\author{
\small {\bf Mourad Bellassoued${}^{1,2,3}$ and
Masahiro Yamamoto${}^4$ }\\
\small ${}^1$ University of Carthage,
Faculty of Sciences of Bizerte,\\
\small Department of Mathematics, 7021 Jarzouna Bizerte, Tunisia.\\
 \small ${}^2$ F\'ed\'eration Denis-Poisson,\\
\small  FR CNRS 2964, B.P. 6759, 45067 Orl\'eans cedex 2, France.\\
\small ${}^3$ LE STUDIUM \circledR, Institute for Advanced Studies, Orl\'eans, France\\
\small mourad.bellassoued@fsb.rnu.tn \\
\small ${}^4$Department of Mathematical Sciences,\\
\small The University of Tokyo, 3-8-1 Komaba, Meguro, Tokyo 153,
 Japan.\\
\small myama@ms.u-tokyo.ac.jp }
\date{}

\maketitle
\begin{abstract}
According to Biot's paper in 1956, by using the Lagrangian
equations in classical mechanics, we consider a problem of
the filtration of a liquid in porous elastic-deformation media
whose mechanical behavior is described by the Lam\'e system
coupled with a hyperbolic equation. Assuming the null surface
displacement on the
whole boundary, we discuss an inverse source problem of
determining a body force only by observation
of surface traction on a suitable subdomain along a sufficiently
large time interval.
Our main result is a H\"older stability estimate for the
inverse source problem, which is
proved by a new Carleman estimat for Biot's system.
\end{abstract}
\section{Introduction}
In 1956, Biot \cite{[Biot]} presented a three-dimensional theory
for coupled frame-fluid wave propagation in fluid saturated
porous media, treating the solid
frame and the saturating fluid as two separate co-located
coupled continua. Two second order coupled partial differential
equations were derived from this theory.
More precisely, let us consider an open and bounded domain $\Omega$ of $\R^3$
with ${\cal C}^\infty$ boundary $\Gamma=\p\Omega$, and
let $\nu = \nu(x)$ be the unit outward normal vector to
$\partial\Omega$ at $x$.
Given $T>0$,  Biot's equation is written as:
\begin{equation}\label{1.1}
\begin{array}{lll}
\varrho_{11}\p_t^2\uu^{s}+\varrho_{12}\p_t^2\uu^{f}
-\DDelta \uu^{s}(x,t)-\nabla\para{\mathrm{q}\,\dive \uu^{f}}
=F_1, & \cr
\varrho_{12}\p_t^2\uu^{s}+\varrho_{22}\p_t^2\uu^{f}
-\nabla\para{\mathrm{q}\,\dive\, \uu^{s}}
-\nabla\para{\mathrm{r}\,\dive \uu^{f}}=F_2, &
\textrm{in}\,\, Q:=\Omega \times (-T,T)
\end{array}
\end{equation}
with the boundary condition
\begin{equation}\label{1.2}
\uu^{s}(x,t)=0,\quad \uu^{f}(x,t)\cdot\nu=0, \quad
(x,t)\in\Sigma:=\Gamma\times(-T,T)
\end{equation}
and the initial condition
\begin{equation}\label{1.3}
\para{\uu^{s}(x,0),\,\p_t\uu^{s}(x,0)}=\para{0,\,0},
\quad \para{\uu^{f}(x,0),\,\p_t\uu^{f}(x,0)} = \para{0,0},
\quad x\in\Omega
\end{equation}
where $F= (F_1,F_2)^T$ is an external force with
$F_\ell = (F_\ell^1, F_\ell^2, F_\ell^3)^T$, $\ell=1,2$, and
$\DDelta$ is the elliptic second-order linear differential operator
given by
\begin{align}\label{1.4}
\DDelta &\vvv(x) \equiv \mu\Delta\vvv(x)
+ \para{\mu+ \lambda}\para{\nabla\ddd\vvv(x)}\cr
&\qquad +\para{\dive\vvv(x)}\nabla\lambda(x)
+\para{\nabla\vvv+(\nabla\vvv)^T}\nabla\mu(x),\qquad x\in\Omega.
\end{align}
Throughout this paper, $t$ and $x=(x_1,x_2,x_3)$ denote the time
variable and the spatial variable respectively, and
$\uu^{s} = \para{u_1^{s}, u_2^{s}, u_3^{s}}^T$ and
$\uu^{f} = \para{u_1^{f}, u_2^{f}, u_3^{f}}^T$ denote
respectively the solid frame and fluid phase displacement
vectors at the location $x$ and the time $t$.

Here and henceforth $\cdot^T$ denotes the transpose of
matrices under consideration.
We assume that the Lam\'e parameters $\mu$ and $\lambda$ satisfy
$$
\mu(x) > 0, \quad \lambda(x) + \mu(x) > 0,\quad \forall
x\in\overline{\Omega}.
$$
The function $\mathrm{q}(x)>0$, $x\in\overline{\Omega}$, is the
dilatational coupling
factor between the fluid phase and the solid frame.
The coefficient $\mathrm{r}(x)>0$, $x\in\overline{\Omega}$ is
the bulk modulus of
the fluid phase and $\varrho_{11}(x),\,\varrho_{22}(x)>0$,
$x\in\overline{\Omega}$ are the corrected mass densities for
the solid phase and the fluid phase porosity and
$\varrho_{12}(x)$ is the inertial coupling factor
and see H\"orlin an Peter \cite{[HP]}.
\medskip

We assume that the sources terms are given by
\begin{equation}\label{1.5}
F_\ell(x,t)=p_\ell(x)R_\ell(x,t),\quad \ell=1,2,\quad (x,t)\in Q,
\end{equation}
where $p_\ell\in H^2(\Omega)$ is real-valued and  $R_\ell=(R_\ell^1,R_\ell^2,R_\ell^3)^T$ satisfy
\begin{equation}\label{1.6}
\sum_{j=2}^3\para{\norm{\p^j_t R}^2_{L^\infty(Q)}+\norm{\p^j_t\nabla R}^2_{L^\infty(Q)}}\leq C.
\end{equation}
The main subject of this paper is the inverse problem of
determining $p=(p_1,p_2) \in (H^2(\Omega))^2$ uniquely from
observed data of
the displacement vector $\uu=(\uu^s,\uu^f)$ in a subdomain
$\omega\subset\Omega$. It is an important problem, for example,
in mechanics to determine the source $p$ inside a porous body
from measurements of the slide frame and fluid phase
displacements in $\omega$.
%
%
%
\subsection{Inverse problem}
Let $\omega\subset\Omega$ be an arbitrarily given subdomain
such that $\partial\omega \supset \partial\Omega$, i.e.,
$\omega=\Omega\cap \mathscr{V}$ where $\mathscr{V}$ is a neighborhood
of $\Gamma$ in $\R^3$ and let $R(x,t)=(R_1(x,t),R_2(x,t))$ be appropriately
given.  Then we want to determine $p(x)=(p_1(x),p_2(x))$, $x\in\Omega$,
by measurements $\uu_{|\omega\times (-T,T)}$.
\medskip

Our formulation of the inverse problem requires only a finite number
of observations. As for inverse problems for non-stationary Lam\'e
system by infinitely many boundary observations
(i.e., Dirichlet-to-Neumann map), we refer to Rachele \cite{[R]},
for example.
\medskip

For the formulation with a finite number of observations,
Bukhgeim and Klibanov \cite{[BK]} created a method
based on a Carleman estimate and established the uniqueness for
inverse problems of determining spatially varying coefficients
for scalar partial differential
equations.  See also Bellassoued \cite{[Be5]}, \cite{[Be6]},
Bellassoued and Yamamoto \cite{[BeYa]},
\cite{[BeYa2]}, Benabdallah, Cristofol, Gaitan and Yamamoto
\cite{[BCGY]}, Bukhgeim \cite{[B]}, Bukhgeim,
Cheng, Isakov and Yamamoto \cite{[BCIY]}, Imanuvilov and Yamamoto
\cite{[IY1]} - \cite{[IY3]}, Isakov \cite{[I1]}, \cite{[I2]},
Kha\u\i darov \cite{[KH1]}, Klibanov \cite{[K1]},
\cite{[KL]}, Klibanov and Timonov \cite{[KT]},
Klibanov and Yamamoto \cite{[KY]}, Yamamoto \cite{[Y]}.
In particular, as for inverse problems for the isotropic Lam\'e
system, we can refer to Ikehata, Nakamura and Yamamoto
\cite{[INY]}, Imanuvilov, Isakov and Yamamoto \cite{[IY4]},
Imanuvilov and Yamamoto \cite{[IY5]} - \cite{[IY6]},
Isakov \cite{[I1]}, Isakov and Kim \cite{[IK1]}.
\medskip

A Carleman estimate is an inequality for a solution to
a partial differential equation with weighted $L^2$-norm and
effectively yields the unique continuation for
a partial differential equation with non-analytic coefficients.
As a pioneering work concerning a Carleman estimate, we refer
to Carleman's paper \cite{[Carleman]} where what is called a Carleman
estimate was proved and applied it
for proving the uniqueness in the
Cauchy problem for a two-dimensional elliptic equation.
Since \cite{[Carleman]}, the theory of Carleman estimates
has been developed and we refer, for example,
to H\"ormander \cite{[H]} and Isakov \cite{[I2]} for Carleman
estimates for functions
having compact supports (that is, they and their derivatives of suitable
orders vanish on the boundary of a domain).
For Carleman estimates for functions without compact
supports, we refer to Bellassoued and Yamamoto \cite{[BeYa3]},
Fursikov and Imanuvilov
\cite{[FI]}, \cite{[Ima2]},
Lavrent'ev, Romanov and Shishat$\cdot$ski\u\i\,
\cite{[La]}, Tataru \cite{[T]}. Moreover Carleman estimates
have been applied for estimating the energy and see e.g.,
Imanuvilov and Yamamoto \cite{[IY6]}, Kazemi and Klibanov
\cite{[KK]}, Klibanov and Malinsky \cite{[KM]},
Klibanov and Timonov \cite{[KT]}.
%
%
%
%
\subsection{Notations and statement of main results}
In order to formulate our results, we need to
introduce some notations.
For $x_0\in\R^3\backslash\overline{\Omega}$, we define
the following set of the scalar coefficients
\begin{equation}\label{1.7}
\mathscr{C}(m,\theta)
=\set{c\in\mathcal{C}^2(\overline{\Omega}),\,
c(x)>c^*>0,x\in\overline{\Omega},
\,\norm{c}_{\mathcal{C}^2(\overline{\Omega})}\leq m,\,
 \frac{\nabla c\cdot(x-x_0)}{2c}\leq 1-\theta},
 \end{equation}
where the constants $m>0$ and $\theta\in (0,1)$  are given.
\medskip

\noindent \textbf{Assumption A.1}\\
Throughout this paper, we assume that the coefficients
$\para{\varrho_{ij}}_{1\leq i,j\leq 2}$, $\mu$, $\lambda$, $\mathrm{q}$,
$\mathrm{r}
\in \mathcal{C}^2(\overline{\Omega})$ satisfy the following conditions
\begin{eqnarray}\label{1.8}
&&\varrho(x)=\varrho_{11}(x)\varrho_{22}(x)
-\varrho_{12}^2(x)>0,\quad \forall x\in\overline{\Omega},\cr
\cr
&& \lambda(x)\mathrm{r}(x)-\mathrm{q}^2(x)>0,\quad
\forall x\in\overline{\Omega}.
\end{eqnarray}
Let $A(x)=\para{a_{ij}(x)}_{1\leq i,j\leq 2}$ be
the $2\times 2$-matrix given by
\begin{equation}\label{1.9}
A(x)=\frac{1}{\varrho}
\left(
\begin{array}{cc}
\varrho_{22} & -\varrho_{12} \\
-\varrho_{12} & \varrho_{11} \\
\end{array}
\right)
\left(
\begin{array}{cc}
2\mu+\lambda & \mathrm{q} \\
 \mathrm{q} & \mathrm{r} \\
 \end{array}
 \right):=\left(
  \begin{array}{cc}
  a_{11} & a_{12} \\
  a_{21} & a_{22} \\
  \end{array}
  \right).
\end{equation}
By (\ref{1.8}), we can prove that $\para{a_{ij}(x)}
_{1\leq i,j\leq 2}$ is a positive definite matrix
on $\overline{\Omega}$.
\medskip

\noindent \textbf{Assumption A.2:}\\
Let $A(x)$ have two distinct positive eigenvalues:
$\mu_2(x),\,\mu_3(x)>0$, $\mu_2(x)\ne \mu_3(x)$.
Moreover setting, $\mu_1:=\para{\varrho^{-1}\varrho_{22}}\mu$,
we assume
\begin{equation}\label{1.10}
\mu_1, \mu_2,\,\mu_3\in \mathscr{C}(m,\theta).
\end{equation}
\medskip

\noindent \textbf{Assumption A.3:}\\
We assume that the solution
$\uu=(\uu^s,\uu^f)$ satisfies the
a priori boundedeness and regularity:
\begin{equation}\label{1.11}
\uu\in H^5(Q),\quad \norm{u}_{H^5(Q)}\leq M_0,
\end{equation}
for some positive constant $M_0$.
\medskip

Before stating the main result on the stability for the
inverse source problem, we present Theorem 1.1
on the unique existence of strong solution to
(\ref{1.1})-(\ref{1.2}) with initial condition:
$$
\uu(\cdot,0) = \uu_0\quad \textrm{and}\quad \partial_t\uu(\cdot,0) = \uu_1.
$$
Let $V(\Omega)=(H^1(\Omega))^3\times H(\dive,\Omega)$, where
\begin{equation}\label{1.12}
H(\dive,\Omega)=\set{\uu\in (L^2(\Omega))^3;\quad\dive\,
\uu \in L^2(\Omega)}.
\end{equation}
The norm in $V(\Omega)$ is chosen as follows
$$
\Vert (\vvv^1, \vvv^2)\Vert_{V(\Omega)}^2
= \norm{\vvv^1}^2_{H^1(\Omega)}
+\norm{\vvv^2}^2_{L^2(\Omega)}+\norm{\dive\vvv^2}^2_{L^2(\Omega)},
\quad \vvv = (\vvv^1, \vvv^2) \in V(\Omega).
$$
\begin{Theo}\label{T0}
Let $F\in H^1(-T,T;L^2(\Omega))$, $(\uu_0,\uu_1)\in
(H^2(\Omega)\cap H^1_0(\Omega))^6\times
(H^1(\Omega))^6$. Then there exists a unique solution $\uu(x,t)
=\para{\uu^s(x,t),\uu^f(x,t)}$ of (\ref{1.1})-(\ref{1.2}) with
initial data $(\uu_0,\uu_1)$ such that
\begin{eqnarray}\label{1.13}
&&\uu^s\in\mathcal{C}([-T,T];H^2(\Omega)\cap H^1_0(\Omega))
\cap\mathcal{C}^1([-T,T];H^1(\Omega))
\cap \mathcal{C}^2([-T,T];L^2(\Omega))\cr
&&\uu^f\in \mathcal{C}^2([-T,T];L^2(\Omega)),
\quad\dive \uu^f\in\mathcal{C}([-T,T];H^1(\Omega))
\cap\mathcal{C}^1([-T,T];L^2(\Omega)).
\end{eqnarray}
In particular there exists a constant $C>0$ such that
$$
\Vert \uu\Vert_{\mathcal{C}^2([-T,T];L^2(\Omega))}
\le C(\Vert F\Vert_{H^1(-T,T;L^2(\Omega))}
+ \Vert \uu_0\Vert_{H^2(\Omega)}
+ \Vert \uu_1\Vert_{H^1(\Omega)}).
$$
Moreover, if $F=0$, then the energy of the solution
$\uu=(\uu^s,\uu^f)$ given by
$$
E(t)= \frac{1}{2}\int_\Omega \para{M(x)\p_t \uu\cdot\p_t \uu
+\lambda|\dive \uu^s|^2+2\mu|\varepsilon( \uu^s)|^2
+\mathrm{r}|\dive \uu^f|^2+2\mathrm{q}(\dive \uu^f)(\dive \uu^s)}dx
$$
is conserved, that is,
$$
E(t)=E(0),\quad \forall t\geq 0.
$$
Here $M(x)=\para{\varrho_{ij}(x)I_3}_{1,\leq i,j\leq 2}$ and
$\varepsilon(\vvv)=\frac{1}{2} \para{\nabla\vvv+(\nabla\vvv)^T}$.
\end{Theo}

The proof is based on the Galerkin
method and see Santos \cite{[S]} for the case $n=2$.
For completeness we will give a proof for dimension $3$
in Section 4.\\

In order to formulate our stability estimates for the inverse
problem we introduce some notations.\\
Let $\vartheta:\overline{\Omega}\To\R$ be
the strictly convex function given by
\begin{equation}\label{1.14}
\vartheta(x)=\abs{x-x_0}^2,\quad x\in\overline{\Omega}.
\end{equation}
Set
\begin{equation}\label{1.15}
\displaystyle D^2=\max_{x\in\overline{\Omega}}\vartheta(x),\quad
\displaystyle d^2=\min_{x\in\overline{\Omega}}\vartheta(x),\quad
\displaystyle D_0^2=D^2-d^2.
\end{equation}
By Assumption A.2, there exist constants $c_j^* > 0$, $j=1,2,3$ such that
$\mu_j(x)>c^*_j>0$ for all $x\in\overline{\Omega}$,
$j=1,2,3$.  Let $c_0^*=\min\{c_1^*,c_2^*,c_3^*\}$.
We choose $\beta > 0$ such that
\begin{equation}\label{1.16}
\beta + \frac{mD_0}{\sqrt{c_0^*}}\sqrt{\beta} < \theta c_0^*,
\qquad
c_0^*d^2 - \beta D^2 > 0.
\end{equation}
Here we note that since $x_0 \not\in \overline{\Omega}$, such
$\beta>0$ exists.\\
We set
\begin{equation}\label{1.17}
T_0=\frac{D_0}{\sqrt{\beta}}.
\end{equation}
The main results of this paper can be stated as follows:
\begin{Theo}(Stability)\label{T1}
Assume (A.1), (A.2), and (A.3). Let $T > T_0$ and
$\uu$ be the solution of
(\ref{1.1})-(\ref{1.2}) and (\ref{1.3}).
Moreover let assume that $\Phi_j(x):=R_j(x,0)$ satisfy
\begin{equation}\label{1.18}
\Phi_j(x)\cdot (x-x_0)\neq 0\quad \textrm{ for all }\,
x\in\overline{\Omega}.
\end{equation}
Let $M>0$. Then there exist constants $C>0$ and $\kappa\in (0,1)$
such that the following estimate holds:
\begin{equation}\label{1.19}
\norm{p_1}^2_{H^1_0(\Omega)}+\norm{p_2}^2_{H^1_0(\Omega)}
\leq C \mathcal{E}_\omega(\uu)^\kappa
\end{equation}
for any $p_\ell\in H^2(\Omega)$, $\ell=1,2$, such that
$\norm{p_\ell}_{H^2}\leq M$ and $p_\ell=0$, $\nabla p_\ell=0$ on $\Gamma$.
Here
$$
\mathcal{E}_\omega(\uu) = \sum_{j=2}^3\norm{\p_t^{j}\uu}^2
_{H^2(\omega\times(-T,T))}.
$$
\end{Theo}
By Theorem \ref{T1}, we can readily derive
the uniqueness in the inverse problem:
\begin{Corol}
Under the assumptions in Theorem \ref{T1}, we have the
uniqueness:\\
Let $\uu=(\uu^s,\uu^f)$ satisfy Biot's
system (\ref{1.1})-(\ref{1.3}) such that
$\uu(x,t)=0,\, (x,t)\in\omega\times (-T,T)$. Then $p_1(x)=p_2(x)
= 0$ for all $x\in\Omega$ and $\uu(x,t)=0$ in $Q$.
\end{Corol}
The remainder of the paper is organized as follows.
In section 2, we give a Carleman estimate for the
Biot's system. In section 3 we prove
Theorem \ref{T1}. Section 4 is devoted
to the proof of Theorem \ref{T0}.
\section{Carleman estimate for Biot's system}
\setcounter{equation}{0}
In this section we will prove a Carleman estimate for
Biot's system, which is interesting of itself.
In order to formulate our Carleman estimate, we
introduce some notations. Let $\vartheta:\overline{\Omega}\To\R$
be the strictly convex function given by (\ref{1.14}),
where $x_0\notin\overline{\Omega}$.

We define two functions $\psi,\varphi:\Omega\times\R\To\R$
of class $\mathcal{C}^\infty$ by
\begin{equation}\label{2.1}
\begin{array}{ll}
\psi(x,t)=\abs{x-x_0}^2-\beta\abs{t}^2 \quad
\textrm{for all } x\in \Omega,\quad -T\leq t\leq T,\cr
\\
\varphi(x,t)=e^{\gamma\psi(x,t)},\quad \gamma>0,
\end{array}
\end{equation}
where $\displaystyle T>T_0$. Therefore, by (\ref{1.17})
and (\ref{1.15}), we have
\begin{equation}\label{2.2}
\varphi(x,0)\geq d_0,\quad \varphi(x,\pm T)< d_0
\end{equation}
with $d_0=\exp(\gamma d^2)$. Thus, for given $\eta>0$, we can
choose sufficiently small $\ep=\ep(\eta)$ such that
\begin{equation}\label{2.3}
\varphi(x,t)\leq d_0 - \eta\equiv d_1  \quad
\textrm{for all }\,\,
(x,t)\in \set{(x,t)\in Q;\,\,\abs{t}>T-2\ep},
\end{equation}
$$
\varphi(x,t) \ge d_0 - \frac{\eta}{2} \equiv d_2
\quad \textrm{for all }\,\,
(x,t)\in \set{(x,t)\in Q;\,\,\abs{t}<\ep}.
$$
Let $(\uu^s,\uu^f)$ satisfy Biot's system
\begin{equation}\label{2.4}
\begin{array}{lll}
\varrho_{11}\partial_t^2\uu^{s}(x,t)
+ \varrho_{12}\partial_t^2\uu^{f}
-\DDelta \uu^{s}(x,t)-\nabla\para{\mathrm{q}\,\dive \uu^{f}}
=F_1, & \cr
\varrho_{12}\partial_t^2\uu^{s}(x,t)
+ \varrho_{22}\partial_t^2\uu^{f}
-\nabla\para{\mathrm{r}\,\dive\,
\uu^{s}}-\nabla\para{\mathrm{q}\,\dive \uu^{f}}=F_2, &
\textrm{in}\,\, Q.
\end{array}
\end{equation}
The following theorem is a Carleman estimate for Biot's
system (\ref{2.4}).
\begin{Theo}\label{T2}
There exist $\tau_* > 0$ and $C>0$  such that
the following estimate holds:
\begin{multline}\label{2.5}
\int_Q\tau\para{|\nabla_{x,t}\uu^s|^2
+|\nabla_{x,t}(\dive \uu^s)|^2
+|\nabla_{x,t}(\dive \uu^f)|^2}e^{2\tau\varphi}dxdt\cr
+ \int_Q\tau^3\para{|\uu^s|^2+|\dive \uu^s|^2
+|\dive \uu^f|^2}e^{2\tau\varphi}dxdt
\leq C\int_Q\para{|F|^2+|\nabla F|^2}e^{2\tau\varphi}dxdt
\end{multline}
for any $\tau\geq\tau_*$ and any solution $(\uu^s,\uu^f)
\in (H^2(Q))^6$ to (\ref{2.4}) which is supported in a
fixed compact set
$K\subset \textrm{int}(Q)$.
\end{Theo}
In order to prove Theorem \ref{T2}, we use a Carleman estimate
for a coupling hyperbolic system, which we discuss in the next
subsection.
\subsection{Carleman estimate for a hyperbolic system}
First we recall the following Carleman estimate for a
scalar hyperbolic equation.  As for the proof, we refer to
Bellassoued and Yamamoto \cite{[BeYa3]},
and Imanuvilov and Yamamoto \cite{[IY3]} for example.
\begin{Lemm}\label{L2.1}
Let $c\in \mathscr{C}(m,\theta)$. There exist constants $C>0$
and $\tau_*>0$ such that the following Carleman estimate holds:
$$
C\int_{Q}\!e^{2\tau\varphi}\para{\tau\abs{\nabla_{x,t} y}^2
+\tau^3\abs{y}^2}\! dx dt
\leq  \int_Q\! e^{2\tau\varphi}\abs{(\p_t^2-
c\Delta)y}^2 dx dt
$$
whenever $y\in H^2(Q)$ is supported in a fixed
compact set $K\subset\mathrm{ int}(Q)$ and any $\tau\geq\tau_*$.
\end{Lemm}

Let $v=(v_1,v_2)\in (H^2(\Omega))^2$ satisfy the following hyperbolic
system
\begin{equation}\label{2.6}
\left\{
\begin{array}{lll}
\p_t^2v_{1}- b_{11}(x)\Delta v_{1}
- b_{12}(x)\Delta v_{2}=g_1 & \textrm{in}\, Q\cr
\p_t^2v_{2}- b_{21}(x)\Delta v_{1}- b_{22}(x)\Delta v_{2}
=g_2 & \textrm{in}\, Q,
\end{array}
\right.
\end{equation}
for $g=(g_1,g_2)\in (L^2(Q))^2$.
We assume that the matrix $B(x) = \para{b_{ij}(x)}
_{1\leq i,j\leq 2}$
has two distinct positive eigenvalues $c_1, c_2\in
\mathscr{C}(m,\theta)$.  Then, by Lemma 2.1, we have the
following Carleman estimate.
\begin{Lemm}\label{L2.2}
There exist constants $C>0$ and $\tau_*>0$ such that
the following Carleman estimate holds
$$
C\int_{Q}\!e^{2\tau\varphi}\para{\tau\abs{\nabla_{x,t} v}^2
+\tau^3\abs{v}^2}\! dx dt\leq  \int_Q\! e^{2\tau\varphi}\abs{g}^2
dx dt
$$
for any $\tau\geq\tau_*$, whenever $v\in H^2(Q)$ is a solution of
(\ref{2.6}) and supported in a fixed compact set $K\subset \textrm{int}(Q)$.
\end{Lemm}
\begin{Demo}{}
The system (\ref{2.6}) can be written in the equivalent form
\begin{equation}\label{2.7}
\p_t^2 v- B(x)\Delta v=g \quad \textrm{in}\,\,Q.
\end{equation}

By the assumption on $B(x)$, there exists a matrix
$P(x)$ such that
$$
\para{P^{-1}BP}(x)=\textrm{Diag}\para{c_1(x),c_2(x)}
= \Lambda (x),\quad x\in\Omega.
$$
Therefore system (\ref{2.7}) can be written in an equivalent
form:
$$
\p_t^2\widetilde{v}-\Lambda(x)\Delta \widetilde{v}
= \widetilde{g}+\mathscr{B}_1(x,\p)v,
$$
where
\begin{equation}\label{2.8}
\widetilde{v}(x,t)=P^{-1}(x)v(x,t),\quad
\widetilde{g}(x,t)=P^{-1}(x)g(x,t),
\end{equation}
and $\mathscr{B}_1$ is a first-order differential operator.\\
Since $c_j\in \mathscr{C}(m,\theta)$ for $j=1,2$, we can apply
Lemma 2.1 for the two components of $\widetilde{v}$ and obtain
$$
C\int_{Q}\!e^{2\tau\varphi}\para{\tau\abs{\nabla_{x,t}
\widetilde{v}}^2+\tau^3\abs{\widetilde{v}}^2}\! dx dt
\leq \int_Q\! e^{2\tau\varphi}|\widetilde{g}|^2 dx dt
+ \int_Q\! e^{2\tau\varphi}\para{|v|^2+|\nabla v|^2}dxdt
$$
and, by (\ref{2.8}), we easily obtain
$$
\abs{v(x,t)}\leq C\abs{\widetilde{v}(x,t)},\quad
\abs{\nabla v(x,t)}\leq C\para{\abs{\nabla\widetilde{v}(x,t)}
+\abs{\widetilde{v}(x,t)}},\quad
|\widetilde{g}(x,t)|\leq C\abs{g(x,t)}
$$
for $(x,t) \in Q$.
This completes the proof.
\end{Demo}
\subsection{Proof of the Carleman estimate for Biot's system}
In this section, we derive a global Carleman estimate for
a solutions of system (\ref{2.4}).
We consider the $6\times 6$-matrix
\begin{equation}\label{2.9}
M(x)=\left(
              \begin{array}{cc}
                \varrho_{11}(x)I_3 & \varrho_{12}(x)I_3 \\
                \varrho_{12}(x)I_3 & \varrho_{22}(x)I_3 \\
              \end{array}
            \right).
\end{equation}
Here $I_3$ is the $3\times 3$ identity matrix.  Then by
Assumption A.1, we have
$$
M^{-1}(x) = \frac{1}{\varrho}\left(
              \begin{array}{cc}
                \varrho_{22}(x)I_3 & -\varrho_{12}(x)I_3 \\
                -\varrho_{12}(x)I_3 & \varrho_{11}(x)I_3 \\
              \end{array}
            \right).
$$

Let $v^s=\dive \uu^s$, $v^f=\dive \uu^f$, $v = (v^s, v^f)$
and $\www^s=\curl \uu^s$.
Put $G=M^{-1}F$, $F=(F_1,F_2)^T$ and apply $M^{-1}$ to
system (\ref{2.4}), we obtain
\begin{equation}\label{2.10}  \begin{array}{lll}
\p_t^2\uu^s -   \mu_1\Delta \uu^s-\para{\mu_1+\lambda_1}
\nabla\para{\dive \uu^s}-\mathrm{q}_1\nabla\para{\dive \uu^f}
= G_1+\mathscr{R}_1\uu^s+\mathscr{R}_0\dive \uu^f,
&\textrm{in}\,Q\cr
\p_t^2\uu^f +\mu_2\Delta \uu^s-\mathrm{r}_2\nabla\para{\dive \uu^f}
-\mathrm{q}_2\nabla\para{\dive \uu^s}
= G_2+\mathscr{R}'_1\uu^s+\mathscr{R}'_0\dive \uu^f,
& \textrm{in}\,Q.
  \end{array}
\end{equation}
Here, $\mathscr{R}_j,\,\mathscr{R}'_j$, $j=0,1$ are
differential operators of order $j$ with coefficients
in $L^{\infty}(Q)$, and
\begin{align}\label{2.11}
&& \mu_1=\varrho^{-1}\mu\varrho_{22},\quad
\lambda_1=\varrho^{-1}\para{\lambda\varrho_{22}-\mathrm{q}\varrho_{12}},
\quad
q_1=\varrho^{-1}\para{\mathrm{q}\varrho_{22}-\mathrm{r}\varrho_{12}}\cr
&&\mu_2=\varrho^{-1}\mu\varrho_{12} ,\quad
\mathrm{q}_2=\varrho^{-1}\para{\mathrm{q}\varrho_{11}-(\mu+\lambda)
\varrho_{12}},
\quad \mathrm{r}_2=\varrho^{-1}\para{\mathrm{r}\varrho_{11}
-\mathrm{q}\varrho_{12}}.
\end{align}
Henceforth $\mathscr{P}_j$, $j=1,...,4$ denote some
first-order operators with $L^{\infty}(Q)$-coefficients.\\
We apply $\dive $ to the equations in (\ref{2.10}),
and can derive the following two equations:
\begin{align}\label{2.12}
\p_t^2v^s-a_{11}\Delta v^s-a_{12}\Delta v^f
&=\dive G_1+\mathscr{P}_1(v^f,v^s,\uu^s,\www^s)\cr
\p_t^2v^f-a_{21}\Delta v^s-a_{22}\Delta v^f
&=\dive G_2+\mathscr{P}_2(v^f,v^s,\uu^s,\www^s),
\end{align}
where $\para{a_{ij}}_{1\leq i,j\leq 2}$ is given by
(\ref{1.9}).
We apply the $\curl$ to the first equation (\ref{2.10}) to
obtain
\begin{equation}\label{2.13}
\p_t^2\www^s-\mu_1\Delta \www^s=\curl G_1
+\mathscr{P}_3(v^f,v^s,\uu^s,\www^s)
\end{equation}
and
\begin{equation}\label{2.14}
\p_t^2\uu^s-\mu_1\Delta \uu^s= G_1
+\mathscr{P}_4(v^f,v^s,\uu^s,\www^s).
\end{equation}
Applying Lemma \ref{L2.2} to system (\ref{2.12}), we have
for $v=(v^s,v^f)$
\begin{multline*}
C\int_{Q}\!e^{2\tau\varphi}\para{\tau\abs{\nabla_{x,t} v}^2
+\tau^3\abs{v}^2}\! dx dt
\leq \int_Q\! e^{2\tau\varphi}\para{\abs{F}^2+\abs{\nabla F}^2}
dxdt\cr
+\int_Qe^{2\tau\varphi}\para{\abs{\uu^s}^2+\abs{\www^s}^2
+\abs{\nabla \uu^s}^2+\abs{\nabla \www^s}^2}dxdt.
\end{multline*}
Applying Lemma \ref{L2.1} to (\ref{2.13}) and (\ref{2.14}), we obtain
\begin{multline*}
C\int_{Q}\!e^{2\tau\varphi}\para{\tau\abs{\nabla_{x,t} \www^s}^2
+\tau^3\abs{\www^s}^2+ \tau\abs{\nabla_{x,t} \uu^s}^2
+\tau^3\abs{\uu^s}^2}\! dx dt\cr
\leq\int_Q\! e^{2\tau\varphi}\para{\abs{F}^2+\abs{\nabla F}^2}
dxdt
+\int_Qe^{2\tau\varphi}\para{\abs{v}^2
+\abs{\nabla v}^2}dxdt.
\end{multline*}
Therefore, for $\tau$ sufficiently
large, we obtain (\ref{2.5}).
This completes the proof of Theorem \ref{T2}.
\section{Proof of Theorem 1.2}
\setcounter{equation}{0}
In this section we prove the stability (Theorem 1.2) for the
inverse source problem.
\medskip

For the proof, we apply the method in Imanuvilov and Yamamoto
\cite{[IY2]} which modified the argument in \cite{[BK]} and proved
the stability for an inverse coefficient problem for a hyperbolic equation.
For it, the Carleman (Theorem 2.1) is a key.
\subsection{Modified Carleman estimate for Biot's system}
Let $\omega_T = \omega \times (-T,T)$.
We modify Theorem 2.1 for functions which vanish at
$\pm T$ with first $t$-derivatives.
\begin{Lemm}\label{L3.1}
There exist positive constants $\tau_*$, $C>0$ and $C_0>0$
such that the following inequality holds:
\begin{multline}\label{3.1}
\int_Q\tau\para{|\nabla_{x,t}\vvv^s|^2
+|\nabla_{x,t}(\dive \vvv^s)|^2
+|\nabla_{x,t}(\dive \vvv^f)|^2}e^{2\tau\varphi}dxdt\cr
+\int_Q\tau^3\para{|\vvv^s|^2
+|\dive \vvv^s|^2+|\dive \vvv^f|^2}e^{2\tau\varphi}dxdt
\leq C\int_Q\para{|G|^2+|\nabla G|^2}e^{2\tau\varphi}dxdt\cr
+ Ce^{C_0\tau}\norm{\vvv}^2_{H^2(\omega_T)}
\end{multline}
for any $\tau\geq \tau_*$ and any $\vvv=(\vvv^s,\vvv^f)
\in H^2(Q)$ satisfying, for $G=(G_1,G_2)$
\begin{equation}\label{3.2}
\begin{array}{lll}
\varrho_{11}\p_t^2\vvv^{s}+\varrho_{12}\p_t^2\vvv^{f}
-\DDelta \vvv^{s}
-\nabla\para{\mathrm{q}\,\dive \vvv^{f}}
&=G_1,\cr
\varrho_{12}\p_t^2\vvv^{s}+\varrho_{22}\p_t^2\vvv^f
-\nabla\para{\mathrm{q}\,\dive\, \vvv^{s}}
-\nabla\para{\mathrm{r}\,\dive \vvv^{f}}&=G_2
\quad \mbox{in $Q$}
\end{array}
\end{equation}
such that
\begin{equation}\label{3.3}
\p^j_t\vvv(x,\pm T)=0 \quad\textrm{for all}
\quad x\in\Omega,\, j=0,1.
\end{equation}
\end{Lemm}
\begin{Demo}{}
Let $\omega^0\subset\omega$. In order to apply Carleman
estimate (\ref{2.5}), we introduce a cut-off function $\xi$
satisfying
$0\leq\xi\leq 1$, $\xi\in {\cal C}^\infty(\R^3)$, $\xi=1$
in $\overline{\Omega\backslash \omega^0}$ and
$\textrm{Supp}\thinspace \xi\subset\Omega$.
Let $\vvv\in H^2(Q)$ satisfy (\ref{3.2}) and (\ref{3.3}). Put
$$
\www(x,t)=\xi(x)\vvv(x,t),\quad (x,t)\in Q,
$$
and let $Q_0=(\Omega\backslash\overline{\omega})\times(-T,T)$.
Noting that $\www\in H^2(Q)$ is compactly supported
in $Q$ and $\www=\vvv$ in $Q_0$ and applying
Carleman estimate (\ref{2.5}) to $\www$, we obtain
\begin{multline*}
\int_{Q_0}\tau\para{|\nabla_{x,t}\vvv^s|^2
+|\nabla_{x,t}(\dive \vvv^s)|^2
+|\nabla_{x,t}(\dive \vvv^f)|^2}e^{2\tau\varphi}dxdt\cr
+ \int_Q\tau^3\para{|\vvv^s|^2+|\dive \vvv^s|^2
+|\dive \vvv^f|^2}e^{2\tau\varphi}dxdt
\leq C\int_Q\para{|G|^2+|\nabla G|^2}e^{2\tau\varphi}dxdt\cr
+ C\int_Q\vert\mathcal{Q}_2\vvv\vert^2
e^{2\tau\varphi}dxdt
\end{multline*}
for any $\tau\geq \tau_*$.  Here $\mathcal{Q}_2$ is
a differential operator of order $2$ whose coefficients
are supported in $\omega$.\\
Therefore
\begin{multline*}
\int_Q\tau\para{|\nabla_{x,t}\vvv^s|^2
+|\nabla_{x,t}(\dive \vvv^s)|^2+|\nabla_{x,t}(\dive \vvv^f)|^2}
e^{2\tau\varphi}dxdt\cr
+ \int_Q\tau^3\para{|\vvv^s|^2+|\dive \vvv^s|^2
+|\dive \vvv^f|^2}e^{2\tau\varphi}dxdt
\leq C\int_Q\para{|G|^2+|\nabla G|^2}e^{2\tau\varphi}dxdt \cr
+ Ce^{C_0\tau}\norm{\vvv}^2_{H^2(\omega_T)}.
\end{multline*}
This completes the proof of the lemma.
\end{Demo}
By $\NN$ we denote
\begin{multline}\label{3.4}
\NN =\int_Q \tau \para{|\nabla_{x,t}\vvv^s|^2
+|\nabla_{x,t}(\dive \vvv^s)|^2
+|\nabla_{x,t}(\dive \vvv^f)|^2}
e^{2\tau\varphi}dxdt\cr
+ \int_Q\tau^3\para{|\vvv^s|^2+|\dive \vvv^s|^2+|\dive \vvv^f|^2}
e^{2\tau\varphi}dxdt
\end{multline}
where $\vvv = (\vvv^s, \vvv^f)$.\\
Now, we recall (\ref{2.2}) and (\ref{2.3}) for the definition of
$d_0$, $\eta$ and $\ep$ and we introduce a cut-off
function $\zeta$ satisfying $0\leq\zeta\leq 1$,
$\zeta\in {\cal C}^\infty(\R)$ and
\begin{equation}\label{3.6}
\zeta=1 \quad \textrm{in $(-T+2\ep,T-2\ep)$}, \quad
\textrm{Supp}\, \zeta\subset(-T+\ep,T-\ep).
\end{equation}
Finally we denote by $\widetilde{\vvv}$ the function
\begin{equation}\label{3.7}
\widetilde{\vvv}(x,t)=\zeta(t)(\vvv^s, \vvv^f)(x,t),
\quad (x,t)\in Q.
\end{equation}
\begin{Lemm}\label{L3.2}
There exist positive constants $\tau_*$, $C$ and $C_0$
such that the following inequality holds:
$$
C\NNN\leq  \int_Q
 \para{
 \abs{ F}^2 + \abs{ \nabla F}^2}
e^{2\tau\varphi}dxdt
+e^{C_0\tau}\norm{\vvv}_{H^2(\omega_T)}^2+
e^{2d_1\tau}\norm{\vvv}^2_{H^1(-T,T;H^1(\Omega))}
$$
for any $\tau\geq \tau_*$ and any $\vvv=(\vvv^s,\vvv^f)\in \para{H^2(Q)}^6$
satisfying
$$\begin{array}{lll}
\varrho_{11}\p_t^2\vvv^{s}+\varrho_{12}\p_t^2\vvv^{f}
-\DDelta \vvv^{s}-\nabla\para{\mathrm{q}\,\dive \vvv^{f}}
=F_1(x,t) & \cr
\varrho_{12}\p_t^2\vvv^{s}+\varrho_{22}\p_t^2\vvv^{f}
-\nabla\para{\mathrm{q}\,\dive\, \vvv^{s}}
-\nabla\para{\mathrm{r}\,\dive \uu^{f}}=F_2(x,t), & \quad (x,t)
\in Q
\end{array}
$$
\end{Lemm}
\begin{Demo}{}
We note that $\widetilde{\vvv}\in \para{H^2(Q)}^6$ and
$$
\begin{array}{lll}
\varrho_{11}\p_t^2\widetilde{\vvv}^{s}
+\varrho_{12}\p_t^2\widetilde{\vvv}^{f}
-\DDelta \widetilde{\vvv}^{s}
- \nabla\para{\mathrm{q}\,\dive \widetilde{\vvv}^{f}}
= \zeta(t) F_1(x,t)+P_1(\vvv,\p_t\vvv),\cr
\varrho_{12}\p_t^2\widetilde{\vvv}^{s}
+\varrho_{22}\p_t^2\widetilde{\vvv}^{f}
-\nabla\para{\mathrm{q}\,\dive\, \widetilde{\vvv}^{s}}
- \nabla\para{\mathrm{r}\,\dive \widetilde{\vvv}^{f}}
=\zeta(t) F_2(x,t)+P_2(\vvv,\p_t\vvv), & \quad (x,t)\in Q,
\end{array}
$$
where $P_1$ and $P_2$ are zeroth-order operators and supported in
$\abs{t}>T-2\ep$.
Therefore, applying Lemma \ref{L3.1} to $\widetilde{\vvv}$
and using (\ref{2.3}), we complete the proof of the lemma.
\end{Demo}
\subsection{Preliminary estimates}

Let $\varphi(x,t)$ be the function defined by (\ref{2.1}). Then
\begin{equation} \label{3.8}
\varphi(x,t)=e^{\gamma\psi(x,t)}=:\rho(x)\alpha(t),
\end{equation}
where $\rho(x)$ and $\alpha(t)$ are defined by
\begin{equation}\label{3.9}
\rho(x)=e^{\gamma\vartheta(x)}\geq d_0,\,\,\forall x\in
\Omega\quad\textrm{and}\quad \alpha(t)=
e^{-\beta\gamma\,t^2}\leq 1,\,\, \forall t\in[-T,T].
\end{equation}
Next we present the following Carleman estimate of
a first-order partial differential operator:
$$
L(x,D)v=\sum_{i=1}^3a_i(x)\p_iv
+ a_0(x)v,\quad x\in\Omega
$$
where
\begin{equation}\label{3.10}
a_0\in{\cal C}(\overline{\Omega}),\quad
a=(a_1,a_2,a_3)\in
\cro{{\cal C}^1(\overline{\Omega})}^3
\end{equation}
and
\begin{equation}\label{3.11}
\abs{a(x)\cdot (x-x_0)}\geq c_0>0,
\quad \textrm{on }\,\overline{\Omega}
\end{equation}
with a constant $c_0>0$.
Then
\begin{Lemm}\label{L3.3}
In addition to (\ref{3.10}) and (\ref{3.11}), we assume
that $\norm{a_0}_{{\cal C}(\overline{\Omega})}\leq M$
and $\norm{a_i}_{{\cal C}^1(\overline{\Omega})}\leq M$,
$1\leq i\leq 3$. Then there exist constants $\tau_*>0$ and
$C>0$ such that
$$
\tau \int_{\Omega}\abs{v(x)}^2e^{2\tau\rho(x)}dx
\leq C\int_\Omega\abs{L(x,D)v(x)}^2e^{2\tau\rho(x)}dx
$$
for all $v\in H^1_0(\Omega)$ and all $\tau>\tau_*$.
\end{Lemm}
The proof is direct by integration by parts and see
e.g., \cite{[IY2]}.
\\

Consider now the following system
\begin{equation}\label{3.12}
\begin{array}{lll}
\varrho_{11}\p_t^2\uu^{s}+\varrho_{12}\p_t^2\uu^{f}
-\DDelta \uu^{s}(x,t)-\nabla\para{\mathrm{q}\,
\dive \uu^{f}} = F_1(x,t), \cr
\varrho_{12}\p_t^2\uu^{s}+\varrho_{22}\p_t^2\uu^{f}
- \nabla\para{\mathrm{q}\,\dive\, \uu^{s}}
- \nabla\para{\mathrm{r}\,\dive \uu^{f}}=F_2(x,t), &
(x,t)\in Q,
\end{array}
\end{equation}
with the boundary condition
\begin{equation}\label{3.13}
\uu^{s}(x,t)=0,\quad \uu^{f}(x,t)\cdot\nu=0, \quad (x,t)\in\Sigma
\end{equation}
and the initial condition
\begin{equation}\label{3.14}
\para{\uu^{s}(x,0),\,\p_t\uu^{s}(x,0)}=(0,0),\quad
\para{\uu^{f}(x,0),\,\p_t\uu^{f}(x,0)}=(0,0), \quad x\in\Omega,
\end{equation}
where the functions $F_1$ and $F_2$ are given by
\begin{equation}\label{3.15}
 F_1(x,t)=p_1(x)R_1(x,t),\quad F_2(x,t)=p_2(x)R_2(x,t).
\end{equation}
We introduce the following notations:
\begin{equation}\label{3.16}
\uu = (\uu^s, \uu^f), \quad
\vvv_j(x,t) = \p_t^j \uu(x,t), \quad (x,t)\in Q,
\thinspace j=0,1,2,3.
\end{equation}
The functions $\vvv_j$, $j=1,2,3$ solve the following system
\begin{equation}\label{3.17}
\begin{array}{lll}
\varrho_{11}\p_t^2\vvv_j^{s}
+\varrho_{12}\p_t^2\vvv_j^{f}-\DDelta \vvv_j^{s}(x,t)
-\nabla\para{\mathrm{q}\,\dive \vvv_j^{f}}
= \p_t^jF_1(x,t),  \cr
\varrho_{12}\p_t^2\vvv_j^{s}+\varrho_{22}\p_t^2\vvv_j^{f}
- \nabla\para{\mathrm{q}\,\dive\, \vvv_j^{s}}
- \nabla\para{\mathrm{r}\,\dive \vvv_j^{f}}
= \p_t^jF_2(x,t), &(x,t)\in Q,
\end{array}
\end{equation}
with the boundary condition
\begin{equation}\label{3.18}
\vvv_j^{s}(x,t)=0,\quad \vvv_j^{f}(x,t)\cdot\nu=0,
\quad (x,t)\in\Sigma.
\end{equation}
We set
$$
\widetilde{\vvv}_j = \zeta \vvv_j,
$$
where $\zeta(t)$ is given by (\ref{3.6}).
We apply Lemma \ref{L3.2} to obtain the following estimate:
\begin{multline}\label{3.19}
C\mathcal{N}_{\tau,\varphi}(\widetilde{\vvv}_j)
\leq \int_Q \para{\vert \p_t^jF\vert^2 + \vert \nabla\p_t^jF\vert^2}
e^{2\tau\varphi}dxdt\cr
+ e^{C_0\tau}\norm{\vvv_j}_{H^2(\omega_T)}^2
+ e^{2d_1\tau}\norm{\vvv_j}^2_{H^1(-T,T;H^1(\Omega))},
\quad j=0,1,2,3,
\end{multline}
provided that $\tau>0$ is large enough.
\begin{Lemm}\label{L3.4}
There exists a positive constant $C>0$ such that the following estimate
$$
\int_\Omega |z(x,0)|^2dx\leq C\int_Q\para{\tau|z(x,t)|^2
+\tau^{-1}|\p_t z(x,t)|^2}dxdt
$$
for any $z\in L^2(Q)$ such that $\p_t z\in L^2(Q)$.
\end{Lemm}
\begin{Demo}{}
Let $\zeta$ be the cut-off function given by (\ref{3.6}).
By direct computations, we have
\begin{eqnarray*}
\int_\Omega\zeta^2(0)|z(x,0)|^2dx
&=&\int_{-T}^{0}\frac{d}{dt}\para{\int_\Omega\zeta^2(t)|z(x,t)|^2
dx}
dt\\
&=&2\int_{-T}^{0}\int_\Omega\zeta^2(t)z(x,t)\p_tz(x,t)dxdt \\
&& +2\int_{-T}^{0}\int_\Omega \zeta'(t)\zeta(t)|z(x,t)|^2dxdt.
\end{eqnarray*}
Then we have
$$
\int_\Omega |z(x,0)|^2dx\leq C\int_Q\para{\tau|z(x,t)|^2
+ \tau^{-1}|\p_tz(x,t)|^2}dxdt.
$$
This completes the proof of the lemma.
\end{Demo}
\begin{Lemm}\label{L3.5}
Let $\phi_\ell(x)=\dive\para{p_\ell(x)\Phi_\ell(x)}$.
Then there exists a constant $C>0$ such that
\begin{multline*}
 \sum_{\ell=1}^2\int_\Omega e^{2\tau\rho}\para{\abs{\phi_\ell(x)}^2
+ \abs{\nabla\phi_\ell(x)}^2}dx\cr
\leq C\para{\mathcal{N}_{\tau,\varphi}(\widetilde{\vvv}_2)
+ \mathcal{N}_{\tau,\varphi}(\widetilde{\vvv}_3)}
+\sum_{\ell=1}^2\int_\Omega\para{\abs{p_\ell}^2+\abs{\nabla p_\ell}^2}
e^{2\tau\rho}dx,
\end{multline*}
provided that $\tau$ is large.
\end{Lemm}
\begin{Demo}{}
We set $\vvv^{(1)} = \vvv_2^s$ and
$\vvv^{(2)} = \vvv_2^f$.
Applying Lemma \ref{L3.4} for
$z_j(x,t)= e^{\tau\varphi(x,t)}\dive \widetilde{\vvv}^{(j)}_2(x,t)$,
$j=1,2$, we obtain the following inequality:
\begin{multline}\label{3.20}
C\tau^2\int_\Omega e^{2\tau\rho}
\sum_{j=1}^2 \abs{\dive\vvv^{(j)}(x,0)}^2dx
\leq \tau^3 \int_Q e^{2\tau\varphi}\sum_{j=1}^2
\abs{\dive\widetilde{\vvv}^{(j)}(x,t)}^2
dxdt\cr
+ \tau \int_Q e^{2\tau\varphi} \sum_{j=1}^2\abs{\p_t\dive
\widetilde{\vvv}^{(j)}(x,t)}^2dxdt\leq \mathcal{N}_{\tau,\varphi}
(\widetilde{\vvv}_2).
\end{multline}
Applying Lemma \ref{L3.4} again with
$w_j(x,t)=e^{\tau\varphi(x,t)}
\nabla\dive(\widetilde{\vvv}^{(j)}(x,t))$, we obtain
\begin{multline}\label{3.21}
C\int_\Omega e^{2\tau\rho}\sum_{j=1}^2
\abs{\nabla\dive\vvv^{(j)}(x,0)}^2dx
\leq\tau \int_Q e^{2\tau\varphi}\sum_{j=1}^2
\abs{\nabla\dive\widetilde{\vvv}^{(j)}(x,t)}^2dxdt\cr
+\tau^{-1} \int_Q e^{2\tau\varphi} \sum_{j=1}^2
\para{\abs{\nabla\dive\partial_t^3\widetilde{\vvv}^s(x,t)}^2
+\abs{\nabla\dive\partial_t^3\widetilde{\vvv}^f(x,t)}^2}
dxdt \cr
\leq \mathcal{N}_{\tau,\varphi}(\widetilde{\vvv}_2)+
\mathcal{N}_{\tau,\varphi}(\widetilde{\vvv}_3).
\end{multline}
Adding (\ref{3.20}) and (\ref{3.21}), we find
\begin{equation}\label{3.22}
\int_\Omega e^{2\tau\rho}\sum_{j=1}^2\para{
\abs{\dive\vvv^{(j)}(x,0)}^2
+ \abs{\nabla\dive\vvv^{(j)}(x,0)}^2}dx
\leq C\para{\mathcal{N}_{\tau,\varphi}(\widetilde{\vvv}_2)
+\mathcal{N}_{\tau,\varphi}(\widetilde{\vvv}_3)}.
\end{equation}
Since
$$
M(x)\para{\vvv^s_2(x,0),\vvv^f_2(x,0)}^T
= \para{p_1(x)\Phi_1(x),p_2(x)\Phi_2(x)}^T, \quad
x \in \overline{\Omega}
$$
 we have
\begin{multline}\label{3.23}
\abs{\phi_\ell(x)}^2+\abs{\nabla\phi_\ell(x)}^2
\leq C\Big(\abs{\vvv^s_2(x,0)}^2+\abs{\nabla\vvv_2^s(x,0)}^2
+ \abs{\vvv^f_2(x,0)}^2 + \abs{\nabla\vvv_2^f(x,0)}^2\cr
+ \vert \dive \vvv^s_2(x,0)\vert^2
+ \vert \nabla \dive\vvv^s_2(x,0)\vert^2
+ \vert \dive \vvv_2^f(x,0)\vert^2
+ \vert \nabla\dive \vvv_2^f(x,0)\vert^2\Big)
\end{multline}
for $x \in \overline{\Omega}$.
On the other hand, using (\ref{3.12}), we obtain
\begin{equation}\label{3.24}
\abs{\vvv^{s}_2(x,0)}^2 + \abs{\nabla\vvv_2^s(x,0)}^2
+ \abs{\vvv^{f}_2(x,0)}^2 + \abs{\nabla\vvv_2^f(x,0)}^2
\le C\sum_{\ell=1}^2\para{\abs{p_\ell}^2+\abs{\nabla p_\ell}^2},
\quad x \in \overline{\Omega}.
\end{equation}
Combining (\ref{3.24}), (\ref{3.23}) and (\ref{3.22}),
we complete the proof of the lemma.
\end{Demo}
\begin{Lemm}\label{L3.6}
There exists a constant $C>0$ such that
$$
\tau\int_\Omega\para{\abs{\nabla p_\ell(x)}^2
+\abs{p_\ell(x)}^2}e^{2\tau\rho}dx
\leq C\int_\Omega \para{\abs{\nabla \phi_\ell(x)}^2
+\abs{\phi_\ell(x)}^2}e^{2\tau\rho(x)}dx
$$
for all large $\tau>0$, $\ell=1,2$.
\end{Lemm}
\begin{Demo}{}
We have
$$
\dive((\p_kp_\ell)(x)\Phi_\ell(x))=\p_k\phi_\ell(x)
-\dive(p_\ell\p_k\Phi_\ell(x))\quad \textrm{for all}\,\, k=1,2,3.
$$
Therefore
\begin{multline}\label{3.25}
\int_\Omega\para{\abs{\dive((\p_k p_\ell)\Phi_\ell)}^2
+\abs{\dive(p_\ell\Phi_\ell)}^2} e^{2\tau\rho}dx
\leq \int_\Omega\para{\abs{\nabla \phi_\ell}^2
+ \abs{\phi_\ell}^2}e^{2\tau\rho(x)}dx\cr
+C\int_\Omega \para{\abs{p_\ell}^2
+ \abs{\nabla p_\ell}^2}e^{2\tau\rho}dx.
\end{multline}
Since $p_\ell=0$ and $\nabla p_\ell=0$ on the boundary $\Gamma$
and $\nabla\Phi_\ell\cdot(x-x_0)\neq 0$, we can apply
Lemma \ref{L3.3} respectively with the choice $v=p_\ell$
and $v=\p_k p_\ell$ to obtain
\begin{equation}\label{3.26}
\tau\int_\Omega\para{\abs{\p_kp_\ell(x)}^2+\abs{p_\ell(x)}^2}
e^{2\tau\rho}dx
\leq C \int_\Omega\para{\abs{\dive((\p_kp_\ell)\Phi_\ell)}^2
+ \abs{\dive(p_\ell\Phi_\ell)}^2} e^{2\tau\rho}dx
\end{equation}
for $\ell=1,2$ and $k=1,2,3$.
Inserting (\ref{3.25}) into the left-hand side of (\ref{3.26})
and choosing $\tau>0$ large, we obtain
$$
\tau\int_\Omega\para{\abs{\nabla p_\ell(x)}^2
+ \abs{p_\ell(x)}^2} e^{2\tau\rho}dx
\leq C\int_\Omega \para{\abs{\nabla \phi_\ell(x)}^2
+ \abs{\phi_\ell(x)}^2}e^{2\tau\rho(x)}dx.
$$
The proof is completed.
\end{Demo}
\subsection{Completion of the proof of Theorem 1.2}
By Lemmata \ref{L3.5} and \ref{L3.6}, we obtain
\begin{multline*}
\tau\sum_{\ell=1}^2\int_\Omega e^{2\tau\rho(x)}
\para{\abs{\nabla p_\ell(x)}^2+\abs{p_\ell(x)}^2}dx
\leq C\sum_{\ell=1}^2\int_\Omega\para{\abs{\nabla p_\ell(x)}^2
+ \abs{p_\ell(x)}^2} e^{2\tau\rho(x)}dx\cr
+ C\para{\mathcal{N}_{\tau,\varphi}(\widetilde{\vvv}_2)
+ \mathcal{N}_{\tau,\varphi}(\widetilde{\vvv}_3)}.
\end{multline*}
Therefore, choosing $\tau>0$ large to absorb the first term on
the right-hand side into the left-hand side and
applying (3.19), we obtain
\begin{multline}\label{3.27}
\tau\sum_{\ell=1}^2\int_\Omega e^{2\tau\rho(x)}
\para{\abs{\nabla p_\ell(x)}^2+\abs{p_\ell(x)}^2}dx
\leq C\sum_{j=2}^3 \biggl(
\int_Q \para{|\p_t^jF\vert^2 + \vert \p_t^j\nabla F\vert^2}
e^{2\tau\varphi}dxdt\cr
+Ce^{C_0\tau}\norm{\vvv_j}_{H^2(\omega_T)}^2
+ Ce^{2d_1\tau}\norm{\vvv_j}^2_{H^1(-T,T;H^1(\Omega))}\biggr)\cr
\leq C\sum_{\ell=1}^2 \int_Q \para{\vert\nabla p_\ell(x)\vert^2
+ \vert p_\ell(x)\vert^2} e^{2\tau\va} dxdt
+ Ce^{C_0\tau} \mathcal{E}_{\omega}(\uu)
+ Ce^{2d_1\tau}M_0.
\end{multline}
Then the first term of the right-hand side of (\ref{3.27}) can be
absorbed into the left-hand side if we take large $\tau>0$.\\
Since $\rho(x)\geq d_0$, we obtain
\begin{equation}\label{3.40}
\sum_{\ell=1}^2\int_\Omega\para{\abs{\nabla p_\ell(x)}^2+\abs{p_\ell(x)}^2}dx
\le C e^{2(d_1-d_0)\tau} + e^{C_0\tau}\mathcal{E}_\omega(\uu)
\leq Ce^{-\epsilon\tau}+e^{C_0\tau}\mathcal{E}_\omega(\uu).
\end{equation}
At the last inequality, we used: By $0 < d_1< d_0$, we can choose $\epsilon > 0$
such that $e^{2(d_1-d_0)\tau} \le e^{-\epsilon\tau}$ for sufficiently
large $\tau > 0$.

\section{Well posedness of the direct problem}
\setcounter{equation}{0}
This section is devoted to the study the existence, uniqueness
and regularity of solutions of the following system:
\begin{equation}\label{4.1}
\begin{array}{lll}
\varrho_{11}\p_t^2\uu^{s}+\varrho_{12}\p_t^2\uu^{f}
-\DDelta \uu^{s}(x,t)-\nabla\para{\mathrm{q}\,\dive \uu^{f}}
=F_1(x,t), & \cr
\varrho_{12}\p_t^2\uu^{s}+\varrho_{22}\p_t^2\uu^{f}
-\nabla\para{\mathrm{q}\,\dive\, \uu^{s}}
-\nabla\para{\mathrm{r}\,\dive \uu^{f}}=F_2(x,t), &
(x,t)\in Q
\end{array}
\end{equation}
with the boundary condition
\begin{equation}\label{4.2}
\uu^{s}(x,t)=0,\quad \uu^{f}(x,t)\cdot\nu=0, \quad
(x,t)\in\Sigma=\Gamma\times(-T,T)
\end{equation}
and the initial condition
\begin{equation}\label{4.3}
\para{\uu^{s}(x,0),\,\uu_t^{s}(x,0)}=\para{\uu_0^s,\,\uu_1^s},
\quad \para{\uu^{f}(x,0),\,\uu_t^{f}(x,0)}=\para{\uu_0^f,\uu_1^f},
\quad x\in\Omega.
\end{equation}
\subsection{Function spaces}
We denote by $\mathscr{D}(\Omega)$ the space of compactly supported,
infinitely differentiable function in $\Omega$ equipped with the
inductive limit topology. We denote by $\mathscr{D}'(\Omega)$ the space
dual to $\mathscr{D}(\Omega)$.  In general, we denote by $X'$ the
space dual to the function space $X$. We denote by $\para{f,g}$ the
inner product in $L^2(\Omega)$ and by $\seq{f,g}$ the value of
$f\in X'$ on $g\in X$. We use usual notations for Sobolev spaces.
If $X$ is a Banach space, then we denote by $L^p(0,T;X)$ the space of
functions $f:(0,T)\To X$ which are measurable, take values in $X$ and
satisfy:
$$
\para{\int_0^T\norm{f(t)}^p_Xdt}^{1/p}=\norm{f}_{L^p(0,T;X)}<\infty
$$
for $1\leq p<\infty$, while
$$
\norm{f}_{L^\infty(0,T;X)}=\textrm{esssup}_{t\in (0,T)}\norm{f(t)}_X<\infty
$$
for $p=\infty$.
It is known that the space $L^p(0,T;X)$ is complete.\\
We define the space
$$
\Ldiv=\set{\uu\in \para{L^2(\Omega)}^3;\,\dive \uu\in L^2(\Omega)},
$$
equipped with the norm
$$
\norm{\uu}_{\Ldiv}=\para{\norm{\uu}_{L^2(\Omega)}^2+\norm{\dive \uu}
_{L^2(\Omega)}^2}^{1/2}.
$$
Let us consider the space
$$
V(\Omega)=\para{H^1(\Omega)^3}\times\Ldiv,
$$
equipped with the norm
$$
\norm{\uu}_{V(\Omega)}=\para{\norm{\uu^2}_{H^1(\Omega)}^2
+ \norm{\uu^2}_{L^2(\Omega)}^2+\norm{\dive \uu^2}_{L^2(\Omega)}^2}^{1/2}.
$$
\subsection{Generalized solution}
We introduce the bilinear form on $V(\Omega)$ by
\begin{multline}\label{4.4}
 B(\uu,\vvv)=\frac{1}{2}\int_\Omega \para{\lambda \dive(\uu^s)\dive(\vvv^s)
+2\mu \para{\varepsilon(\uu^s):\varepsilon(\vvv^s)}
+r\dive(\uu^f)\dive(\vvv^f)}dx\cr
+\frac{1}{2}\int_\Omega q\para{\dive(\uu^f)\dive(\vvv^s)+
\dive(\vvv^f)\dive(\uu^s)}dx
\end{multline}
for any $\uu=(\uu^s,\uu^f)\in V(\Omega)$, $\vvv=(\vvv^s,\vvv^f)\in V(\Omega)$.
We recall that the matrix $M$ is given by (\ref{2.9}).
\begin{defi}
We say that $\uu=(\uu^s,\uu^f)$ is a generalized solution
of problem (\ref{4.1})-(\ref{4.2}), if $\uu\in L^2(0,T;V(\Omega))$
satisfies the initial condition (\ref{4.3}) and the following identity
\begin{equation}\label{4.5}
\para{M \p_t^2\uu(t),\vvv(t)}+B(\uu(t),\vvv(t))=\para{F(t),\vvv(t)},
\quad \mbox{almost all $t\in (0,T)$}
\end{equation}
for any $\vvv\in L^2(0,T;V(\Omega))$.
\end{defi}
We note that in (\ref{4.5}) the integration is only in $x$.

\begin{Lemm}\label{L4.1}
For $\eta>0$, we set
$$
B_\eta(\uu,\vvv)=B(\uu,\vvv)+\eta\para{\uu,\vvv},\qquad \uu,\vvv\in V(\Omega).
$$
Then there exists sufficiently large constant $\eta$ such that
the symmetric bilinear
form $B_\eta$ satisfies
\begin{itemize}
  \item[(i)]
 $\abs{B_\eta(\uu,\vvv)}\leq C_1\norm{\uu}_{V(\Omega)}\norm{\vvv}
_{V(\Omega)}$, for any $\uu,\vvv\in V(\Omega)$,
 \item[(ii)]
$B_\eta(\uu,\uu)\geq C_2\norm{\uu}^2_{V(\Omega)}$, for any $\uu\in V(\Omega)$.
\end{itemize}
\end{Lemm}
\begin{Demo}{}
By (\ref{4.4}) we obtain, for any $\uu,\vvv\in V(\Omega)$
\begin{align}\label{4.6}
\abs{B(\uu,\vvv)}&\leq \para{\norm{\uu^s}_{H^1(\Omega)}
+ \norm{\dive \uu^f}_{L^2(\Omega)}}\para{\norm{\vvv^s}_{H^1(\Omega)}
+ \norm{\dive \vvv^f}_{L^2(\Omega)}}\cr
&\leq C\norm{\uu}_V\norm{\vvv}_V.
\end{align}
Then for any $\eta$, we can derive (i).\\
Now, we note that for a vector $\uu^s\in H^1_0(\Omega)$ we have the
following Korn's inequality
$$
C_1\norm{\uu^s}^2_{H^1(\Omega)}\leq \int_\Omega \varepsilon(\uu^s):
\varepsilon(\uu^s)dx.
$$
Then, for $W=\para{\dive \uu^s,\dive \uu^f}$, we have
$$
B(\uu,\uu)\geq \mu C_1\norm{\uu^s}^2_{H^1(\Omega)}+\frac{1}{2}\int_\Omega M_0W
\cdot W dx
$$
where $M_0$ is the symmetric $2\times 2$-matrix given by
$$
M_0(x)=\left(
  \begin{array}{cc}
    \lambda & q \\
    q & r \\
  \end{array}
\right)\geq \gamma_0 I_2.
$$
Which implies
\begin{align}\label{4.7}
B(\uu,\uu) &\geq \mu C_1\norm{\uu^s}^2_{H^1(\Omega)}+\frac{\gamma_0}{2}
\para{\norm{\dive \uu^s}^2_{L^2(\Omega)}+\norm{\dive \uu^f}^2
_{L^2(\Omega)}}\cr
&\geq \mu C_1\norm{\uu^s}^2_{H^1(\Omega)}+\frac{\gamma_0}{2}
\norm{\dive \uu^f}^2_{L^2(\Omega)}-\frac{\gamma_0}{2}\norm{\uu}^2
_{L^2(\Omega)}\cr
&\geq C_2 \norm{\uu}^2_V-\eta\norm{\uu}^2_{L^2(\Omega)}.
\end{align}
This completes the proof of the lemma.
\end{Demo}
\subsection{Construction of approximate solutions}
Let $\mathscr{A}:\para{L^2(\Omega)}^6\rightarrow \para{L^2(\Omega)}^6$
be the self-adjoint operator defined by
$$
\mathscr{A}u=\left(
     \begin{array}{c}
       \Delta_{\mu,\lambda}\uu^s+\nabla\para{q\dive \uu^f} \\
       \nabla\para{q\dive \uu^s}+\nabla\para{r\dive \uu^f} \\
     \end{array}
   \right).
$$
Then system (\ref{4.1}) can be written as
\begin{equation}\label{4.8}
M\p_t^2\uu-\mathscr{A}\uu=F,\quad (x,t)\in Q
\end{equation}
with initial condition
\begin{equation}\label{4.9}
\uu(x,0)=(\uu_0^s(x),\uu_0^f(x)),\quad \p_t\uu(x,0)=(\uu_1^s(x),\uu_1^f(x))
\end{equation}
and the boundary condition
\begin{equation}\label{4.10}
\uu^s(x,t)=0,\quad \uu^f\cdot\nu=0,\quad (x,t)\in\Sigma.
\end{equation}
Let $(\www_j)_{j\geq 1}$ be a sequence of solutions in $\para{H^2(\Omega)
\cap H^1_0(\Omega)}^6$ such that for all $m\in\N$, $\www_1,...,\www_m$
are linearly independent and all the finite linear combinations of
$(\www_j)_{j\ge 1}$ are dense in $\para{H^2(\Omega)}^6$.\\
We seek approximate  solutions of the problem in the form
\begin{equation}\label{4.11}
\uu_m(t)=\sum_{j=1}^mg_{jm}(t)\www_j.
\end{equation}
The functions $g_{jm}(t)$ are defined by the solution of the system
of ordinary differential equations
\begin{equation}\label{4.12}
\para{M\p_t^2\uu_m,\www_j}+B(\uu_m,\www_j)=\para{F(t),\www_j},\quad 1\leq j\leq m,
\end{equation}
with the initial conditions
\begin{align}\label{4.13}
\uu_m(0)=\uu_{0m} &\rightarrow \uu_0\quad\textrm{in}\quad
\para{H^2(\Omega)\cap H^1_0(\Omega)}^6,\cr
\p_t\uu_m(0)=\uu_{1m} &\rightarrow \uu_1 \quad\textrm{in}\quad
\para{H^1(\Omega)}^6.
\end{align}
The system (\ref{4.12})-(\ref{4.13}) depends on $g_{jm}(t)$ and therefore
has a solution on some segment $[0,t_m]$; see \cite{[LM]}. From a priori
estimates below and the theorem on
continuation of a solution we deduce that it is possible to take $t_m=T$.
\subsection{A priori estimates}
Multiplying (\ref{4.8}) by $g'_{jm}(t)$ and summing over $j$ from 1 to $m$,
we obtain
\begin{equation}\label{4.14}
\para{M\p_t^2\uu_m,\p_t\uu_m}+B(\uu_m,\p_t\uu_m)=\para{F(t),\p_t\uu_m}.
\end{equation}
Hence
\begin{equation}\label{4.15}
\frac{1}{2}\frac{d}{dt}\cro{\norm{M^{1/2}\p_t\uu_m(t)}^2_{L^2(\Omega)}
+B_\eta(\uu_m(t),\uu_m(t))}=(F(t),\p_t\uu_m(t))
+\frac{\eta}{2}\frac{d}{dt}\norm{\uu_m(t)}^2_{L^2(\Omega)}.
\end{equation}
Let
$$
\Phi^2(t)=\norm{M^{1/2}\p_t\uu_m(t)}^2_{L^2(\Omega)}+B_\eta(\uu_m(t),\uu_m(t)).
$$
From (\ref{4.15}) we obtain
\begin{equation}\label{4.16}
\frac{1}{2}\frac{d}{dt}\Phi^2(t)\leq C\cro{\norm{F(t)}^2_{L^2(\Omega)}
+\norm{\p_t\uu_m(t)}^2_{L^2(\Omega)}+\norm{\uu_m(t)}^2_{L^2(\Omega)}}.
\end{equation}
Integrating with respect to $\tau$ from $0$ to $t$, we obtain
\begin{equation}\label{4.17}
\Phi^2(t)\leq C \cro{\norm{F}^2_{L^2(Q)}+\Phi^2(0)
+ \int_0^t\para{\norm{\p_t\uu_m(\tau)}^2_{L^2(\Omega)}
+ \norm{\uu_m(\tau)}^2_{L^2(\Omega)}}}.
\end{equation}
Since
\begin{equation}\label{4.18}
\Phi^2(t)\geq C\para{\norm{\p_t\uu_m(t)}^2_{L^2(\Omega)}
+ \norm{\uu_m(t)}^2_{V(\Omega)}}
\end{equation}
and
\begin{equation}\label{4.19}
\Phi^2(0)\leq C+\norm{\uu_0}^2_{H^2(\Omega)}+\norm{\uu_1}^2_{H^1(\Omega)},
\end{equation}
we have from (\ref{4.18})
\begin{equation}\label{4.20}
\norm{\uu_m(t)}_{V(\Omega)}^2+\norm{\p_t\uu_m(t)}^2_{L^2(\Omega)}\leq R_0+
\int_0^t\para{\norm{\p_t\uu_m(\tau)}^2_{L^2(\Omega)}+ \norm{\uu_m(\tau)}^2
_{L^2(\Omega)}},
\end{equation}
where $R_0=C+\norm{\uu_0}^2_{H^2(\Omega)}
+\norm{\uu_1}^2_{H^1(\Omega)}+\norm{F}^2_{L^2(Q)}$.
By the Gronwall inequality, we conclude that
\begin{equation}\label{4.21}
\norm{\uu_m(t)}_{V(\Omega)}^2+\norm{\p_t\uu_m(t)}^2_{L^2(\Omega)}
\leq R_0
\end{equation}
for all $t\in (0,T)$ and $m\geq 1$.

In order to obtain the second a priori estimate, we observe that
\begin{equation}\label{4.22}
\norm{\p_t^2\uu_m(0)}^2_{L^2(\Omega)}\leq C\cro{\norm{F(0)}^2_{L^2(\Omega)}
+\norm{\uu_m(0)}^2_{H^2(\Omega)}+\norm{\p_t\uu_m(0)}^2_{L^2(\Omega)}}\leq R_1.
\end{equation}
Indeed, multiplying (\ref{4.8}) by $g'_{jm}(0)$, summing over $j$ and
setting $t=0$, we obtain
\begin{equation}\label{4.23}
\para{M\p_t^2\uu_m(0),\p_t^2\uu_m(0)}+B(\uu_m(0),\p_t^2\uu_m(0))
=\para{F(0),\p_t^2\uu_m(0)}.
\end{equation}
Consequently,
\begin{equation}\label{4.24}
\para{M\p_t^2\uu_m(0),\p_t^2\uu_m(0)}=\para{F(0),\p_t^2\uu_m(0)}
+\para{\mathscr{A}\uu_{0m},\p_t^2\uu_m(0)},
\end{equation}
which implies
\begin{equation}\label{4.25}
\norm{\p_t^2\uu_m(0)}^2\leq C\para{\norm{F(0)}^2
_{L^2(\Omega)}+\norm{\uu_{0m}}^2_{H^2(\Omega)}}\leq CR_2.
\end{equation}
Differentiating (\ref{4.14}) with respect to $t$, multiplying by $g_{jm}$
and summing over $j$, we obtain the identity
\begin{equation}\label{4.26}
\frac{1}{2}\frac{d}{dt}\cro{\norm{M^{1/2}\p_t\uu_m(t)}^2_{L^2(\Omega)}
+B_\eta(\uu_m(t),\uu_m(t))}=(F(t),\p^2_t\uu_m(t))
+\frac{\eta}{2}\frac{d}{dt}\norm{\p_t\uu_m(t)}^2_{L^2(\Omega)}.
\end{equation}
Then, we conclude that
\begin{align}\label{4.27}
\norm{\p_t^2\uu_m(t)}^2_{L^2(\Omega)}+\norm{\p_t\uu_m(t)}^2_V & \leq R_2
+\norm{\p_t\uu_m(0)}^2_{L^2(\Omega)}\cr
&+C\int_0^t\para{\norm{\p_t^2\uu_m(\tau)}^2_{L^2(\Omega)}
+\norm{\p_t\uu_m(\tau)}^2_V}d\tau.
\end{align}
By (\ref{4.27}) and the Gronwall inequality, we obtain
\begin{equation}\label{4.28}
\norm{\p_t^2\uu_m(t)}^2_{L^2(\Omega)}+\norm{\p_t\uu_m(t)}^2_V\leq R_1.
\end{equation}
Taking into consideration that $\uu_m=0$ in $\Sigma$, we see
\begin{align}\label{4.29}
\uu_m & \in L^\infty(0,T;V(\Omega)),\quad
\p_t\uu_m  \in L^\infty(0,T;V(\Omega)),\cr
\p_t^2\uu_m & \in L^\infty(0,T;L^2(\Omega)).
\end{align}
\subsection{Passage to the limit}
By (\ref{4.29}), we can extract a sequence from $(\uu_m)_{m\geq 0}$,
which we denote again by $(\uu_m)_m$, such that
\begin{align}\label{4.30}
\uu_m &\rightarrow \uu \quad \textrm{in the weak-star topology
in }\,\, L^\infty(0,T;V(\Omega))\cr
\p_t\uu_m & \rightarrow\p_t\uu\quad  \textrm{in the weak-star topology in }\,\, L^\infty(0,T;V(\Omega))\cr
\p_t^2\uu_m & \rightarrow \p_t^2\uu\quad
\textrm{in the weak-star topology in }\,\, L^\infty(0,T;L^2(\Omega))
\end{align}
and
$$
(\uu_m,\p_t\uu_m)\rightarrow (\uu,\p_t\uu)\quad
\textrm{a.e., on }\,\,\Sigma.
$$
Multiplying (\ref{4.8}) by $\theta\in L^1(0,T)$ and integrating, we have
\begin{equation}\label{4.31}
\int_0^T\para{\para{M\p_t^2\uu_m(t),\www_j}+B(\uu_m,\www_j)}\theta(t)dt
=\int_0^T\para{F(t),\www_j}\theta(t)dt.
\end{equation}
On the other hand
\begin{equation}\label{4.32}
\int_0^TB(\uu_m,\www_j)\theta(t)dt
=-\int_0^T(\uu_m,\mathscr{A}\www_j)\theta(t)dt,
\end{equation}
so that
\begin{equation}\label{4.33}
\lim_{m\to\infty}\int_0^TB(\uu_m,\www_j)\theta(t)dt
= -\int_0^T(\uu,\mathscr{A}\www_j)\theta(t)dt=\int_0^TB(\uu,\www_j)\theta(t)dt.
\end{equation}
Thus, we obtain
\begin{equation}\label{4.34}
\int_0^T\para{\para{M\p_t^2\uu(t),\www_j}+B(\uu,\www_j)}\theta(t)dt
=\int_0^T(F(t),\www_j)\theta(t)dt.
\end{equation}
Taking into account that $\www_j$ are dense in $\para{H^2(\Omega)\cap
H^1_0(\Omega)}^6$ and therefore in $V$, we obtain
\begin{equation}\label{4.35}
\para{M\p_t^2\uu,\vvv}+B(\uu(t),\vvv(t))=(F(t),\vvv),\quad t\in(0,T)
\end{equation}
for all $\vvv\in L^2(0,T;V(\Omega))$.\\
We have $B(\uu(t),\vvv(t))=-(\mathscr{A}\uu(t),\vvv(t))$ for any
$\vvv\in \mathscr{D}(\Omega)$, where the application of the differential
operator $A$ to $\uu$ is in
the distributional sense in $\mathscr{D}'(\Omega)$.  Hence we obtain
\begin{equation}\label{4.36}
M\p_t^2\uu-\mathscr{A}\uu=F,\quad \textrm{in}\,\,\mathscr{D}'(\Omega),
\,\,\textrm{a.e. in} (0,T).
\end{equation}
On the other hand, $\p_t^2\uu$, $\p_t\uu$, $F \in
L^\infty(0,T;L^2(\Omega))$.  Hence (\ref{4.36}) holds in
$L^\infty(0,T;L^2(\Omega))$.\\
The boundary condition (\ref{4.2}) is satisfied
by the choice of the space $V(\Omega)$. We prove that the initial
conditions are satisfied. Suppose
$\theta\in \mathcal{C}^1(0,T)$ and $\theta(T)=0$. For any $j$ we have
\begin{equation}\label{4.37}
\int_0^T\para{\frac{\p}{\p t}(\uu_m-\uu),\www_j}\theta(t)dt
=-\para{\uu_m(0)-\uu(0),\www_j}\theta(0)-\int_0^T\para{\uu_m(t)-\uu(t),\www_j}
\theta'(t)dt.
\end{equation}
Then, by (\ref{4.30}), we have
$$
\lim_{m\to\infty}\abs{\para{\uu_m(0)-\uu(0),\www_j}}=0.
$$
Since $\uu_{0m}(x)=\uu_m(0,x)$ and $\uu_{0m}\rightarrow \uu_0$,
we obtain $\uu(0)=\uu_0$, and can argue similarly for $\uu_1$.\\
Then, we conclude that, there exists a solution $u$ of (\ref{4.1}) such that
\begin{equation}\label{4.38}
\p_t \uu\in L^\infty(0,T;V(\Omega)),
\quad \textrm{and}\quad \p^2_t \uu\in L^\infty(0,T;L^2(\Omega)),
\end{equation}
which implies
\begin{align}\label{4.39}
\uu^s\in \mathcal{C}^1(0,T;H^1_0(\Omega)),&\quad \p_t^2\uu^s
\in \mathcal{C}(0,T;L^2(\Omega))\cr
\uu^f\in \mathcal{C}^1(0,T;H(\dive,\Omega)),&\quad
\p_t^2\uu^f\in \mathcal{C}(0,T;L^2(\Omega)).
\end{align}
On the other hand
\begin{align}\label{4.40}
\nabla\para{\mathrm{q}\dive \uu^s}+\nabla\para{\mathrm{r}\dive \uu^f}
\in\mathcal{C}(0,T;L^2(\Omega)),\cr
\Delta_{\mu,\lambda}\uu^s+\nabla\para{\mathrm{q}\dive \uu^f}
\in \mathcal{C}(0,T;L^2(\Omega)).
\end{align}
Consequently,
$$
\Delta_{\mu,\tild{\lambda}}\uu^s\in\mathcal{C}(0,T;L^2(\Omega)),\quad
\tild{\lambda}=\lambda-\frac{\mathrm{q}^2}{\mathrm{r}}.
$$
Then by the elliptic regularity, $\uu^s\in H^1_0(\Omega)$ yields
$$
\uu^s\in\mathcal{C}(0,T;H^2(\Omega)).
$$
By $\dive \uu^f\in \mathcal{C}(0,T;H^1(\Omega))$, we see
\begin{align}\label{4.41}
\uu^s\in\mathcal{C}(0,T;H^2(\Omega)\cap H^1_0(\Omega))
\cap\mathcal{C}^1(0,T;H^1(\Omega))\cr
\uu^f\in \mathcal{C}^2(0,T;L^2(\Omega)),\quad
\dive \uu^f\in\mathcal{C}(0,T;H^1(\Omega))\cap\mathcal{C}^1(0,T;L^2(\Omega)).
\end{align}
\subsection{Uniqueness}
Let $\uu_1$ and $\uu_2$ be two solutions to (\ref{4.1})-(\ref{4.2}) with
the same initial data, and set $\uu=\uu_1-\uu_2$.
Then for every function $\vvv\in V(\Omega)$, we have
$$
\para{M\p_t^2\uu,\vvv}+B(\uu,\vvv)=0,\quad\forall t\in(0,T).
$$
Since $\p_t\uu\in V(\Omega)$, we may take $\vvv=\p_t\uu$, and this equation
can be reduced to equality
$$
\frac{1}{2}\frac{d}{dt}\cro{\norm{M^{1/2}\p_t\uu}^2_{L^2(\Omega)}
+B_\eta(\uu,\uu)}=\frac{\eta}{2}\frac{d}{dt}\norm{\uu(t)}^2_{L^2(\Omega)}.
$$
Then
$$
\norm{\p_t\uu(t)}^2_{L^2(\Omega)}+\norm{\uu(t)}^2_{V(\Omega)}
\leq C\int_0^t\para{\norm{\p_t\uu(\tau)}^2_{L^2(\Omega)}+\norm{\uu(\tau)}^2
_{L^2(\Omega)}}d\tau.
$$
This implies that $\norm{\uu}_{V(\Omega)}=0=\norm{\p_t\uu}_{L^2(\Omega)}$ and
$\uu_1=\uu_2$ a.e. in $Q$.\\
The proof of Theorem \ref{T0} is completed.

{\bf Acknowledgements:}
Most of this article has been written during the stays of the first author
at Graduate School of
Mathematical Sciences of the University of Tokyo in 2010 and 2011.
The author thanks the school for the hospitality.



\begin{thebibliography}{99}
%
\bibitem{[Be5]}{M.Bellassoued: }
{\it Global logarithmic stability in inverse hyperbolic
problem by arbitrary boundary observation}, Inverse Problems 20 (2004),
1033-1052.

\bibitem{[Be6]}{M. Bellassoued: }{\it Uniqueness and stability in determining
the speed of propagation of second-order hyperbolic equation with variable
coefficients}, Applicable Analysis 83 (2004), 983-1014.

\bibitem{[BeYa]}{M. Bellassoued and M.Yamamoto: }
{\it Logarithmic stability in determination of a coefficient in an acoustic
equation by arbitrary boundary observation}, J. Math. Pures Appl.
85(2006), 193-224.

\bibitem{[BeYa2]}{M. Bellassoued and M.Yamamoto: }
{\it Determination of a coeffcient in the wave
equation with a single measurment}, Applicable
Analysis 87 (2008), 901-920.

\bibitem{[BeYa3]}{M. Bellassoued and M.Yamamoto: }
{\it
Carleman estimate with second large parameter for second
order hyperbolic operators in a Riemannian manifold
and applications in thermoelasticity cases},
Applicable Analysis 91 (2012), 35--67.

\bibitem{[BCGY]}{A. Benabdallah, M. Cristofol,
P. Gaitan and M. Yamamoto:}
{\it Inverse problem for a parabolic system with
two components by measurements of one component}.
Applicable Analysis 88 (2009), 683 -- 709.

\bibitem{[Biot]}{ M.A. Biot: }{\it Theory of propagation
of elastic waves in a fluid-saturated porous solid. I.
Low-frequency range.}
J. Acoust. Soc. Amer. 28 (1956) 168-178, 73.99.

\bibitem{[B]}{A.L.Bukhgeim: }
{\it Introduction to the Theory of Inverse Problems}, VSP,
Utrecht (2000).

\bibitem{[BCIY]}{A.L.Bukhgeim, J.Cheng, V. Isakov
and M.Yamamoto: }
{\it Uniqueness in determining damping corfficients in
hyperbolic equations}, S.Saitoh et al. (eds),
Analytic Extension Formulas and their Applications, 27-46 (2001).

\bibitem{[BK]}{A.L. Bukhgeim and M.V. Klibanov: }
{\it Global uniqueness of class of multidimentional
inverse problems,} Soviet
Math. Dokl. 24 (1981), 244-247.
\bibitem{[Carleman]}{T.Carleman: }
{\it Sur un probl\`eme d'unicit\'e pour les syste\`eme
d'\'equations aux d\'eriv\'ees partielles \`a deux variables
ind\'ependents}, Ark. Mat. Astr. Fys. 2B (1939), 1-9.



\bibitem{[FI]}{A.V. Fursikov and O.Yu. Imanuvilov:}
{\it Controllability of Evolution Equations},
Seoul National University, Seoul (1996).

\bibitem{[HP]}{N.-E. H\"orlin and G.Peter: }{\it Weak, anisotropic symmetric
formulations of Biot's equations for vibro-acoustic modelling of porous
elastic materials}. Internat. J. Numer. Methods Engrg. 84 (2010), no. 12,
1519-1540.

\bibitem{[H]}{L.H\"ormander: }
{\it Linear Partial Differential Operators},
Springer-Verlag, Berlin
(1963).

\bibitem{[INY]}{M. Ikehata, G.Nakamura and M.Yamamoto:}
{\it Uniqueness in inverse problems for the isotropic
Lam\'e system}, J. Math. Sci. Univ. Tokyo 5 (1998), 627-692.

\bibitem{[Ima2]}{O.Yu. Imanuvilov:}
{\it On Carleman estimates for hyperbolic equations},
Asymptotic Analysis 32 (2002), 185--220.

\bibitem{[IIY]}{O.Yu.Imanuvilov, V.Isakov and M.Yamamoto: }
{\it An inverse problem for the dynamical Lam\'e system
with two sets of
boundary data}, Comm. Pure Appl. Math. 56 (2003),
1366-1382.

\bibitem{[IY1]}{O.Yu.Imanuvilov and M.Yamamoto: }
{\it Lipshitz stability in inverse parabolic problems
by Carleman estimate},
Inverse Problems 14 (1998), 1229-1249.

\bibitem{[IY2]}{O.Yu.Imanuvilov and M.Yamamoto: }
{\it Global Lipschitz
stability in an inverse hyperbolic problem by
interior observations},
Inverse Problems 17 (2001), 717-728.

\bibitem{[IY3]}{O.Yu.Imanuvilov and M.Yamamoto: }
{\it Determination of a coefficient in an acoustic
equation with single measurement}, Inverse Problems
19 (2003), 157-171.

\bibitem{[IY4]}{O.Yu.Imanuvilov and M.Yamamoto: }
{\it Carleman estimates for the non-stationary
Lam\'e system and the application to an
inverse problem},
ESIAM:COCV 11 (2005), 1-56.

\bibitem{[IY5]}{O.Yu.Imanuvilov and M.Yamamoto: }
{\it Carleman estimates for the Lam\'e system with
stress boundary condition},
Publ. Research Institute for Mathematical Sciences
43 (2007), 1023--1093.

\bibitem{[IY6]}{O.Yu.Imanuvilov and M.Yamamoto: }
{\it  An inverse problem and an observability
inequality for the Lame system with stress
boundary condition}, Applicable Analysis 88 (2009),
711--733.

\bibitem{[I1]}{V.Isakov: }
{\it A nonhyperbolic Cauchy problem for $\Box_b\Box_c$ and its applications
to elasticity theory},
Comm. Pure and Appl. Math. 39 (1986), 747-767.

\bibitem{[I2]}{V.Isakov: }
{\it Inverse Problems for Partial Differential
Equations}, Springer-Verlag, Berlin (1998, 2005).

\bibitem{[IK1]}{V. Isakov and N. Kim:}
{\it Carleman estimates with second large parameter
and applications to elasticity with residual stress},
Applicationes Mathematicae 35 (2008), 447-465.


\bibitem{[KK]}{M.A.Kazemi and M.V.Klibanov: }
{\it Stability estimates for ill-posed Cauchy
problems involving hyperbolic equations and
inequality}, Applicable Analysis 50 (1993), 93-102.

\bibitem{[KH1]}{A.Kha\u\i darov: }
{\it On stability estimates in multidimentional inverse problems for
differential equation},  Soviet Math. Dokl. 38 (1989), 614-617.

\bibitem{[K1]}{M.V. Klibanov:}
{\it Inverse problems in the "large" and Carleman bounds},
Differential Equations 20 (1984), 755-760.

\bibitem{[KL]}{M.V.Klibanov: }
{\it Inverse problems and Carleman estimates}, Inverse Problems 8 (1992),
575-596.

\bibitem{[KM]}{M.V. Klibanov and J. Malinsky:}
{\it Newton-Kantorovich method for 3-dimensional potential
inverse scattering problem and stability of the hyperbolic
Cauchy problem with time dependent data}, Inverse Problems 7 (1991), 577--595.

\bibitem{[KT]} {M.V. Klibanov and A. Timonov:}
{\it Carleman Estimates for Coefficient Inverse
Problems and Numerical Applications}, VSP, Utrecht (2004).

\bibitem{[KY]}{M.V. Klibanov and M. Yamamoto:}
{\it Lipschitz stability of an inverse problem for
an accoustic equation},
Applicable Analysis 85 (2006), 515-538.

\bibitem{[La]}{M.M.Lavrent'ev, V.G. Romanov and
S.P. Shishat$\cdot$ski\u\i: }
{\it Ill-posed Problems of Mathematics Physics and Analysis},
American Math. Soc., Providence (1986).
%

%
\bibitem{[LM]}{J.-L. Lions and E. Magenes:}
{\it Non-homoeneous Boundary Value Problems
and Applications}, vol. I and II,
Springer-Verlag, Berlin, 1972.
%



\bibitem{[R]}{L.Rachele: }{\it An inverse problem in
elastodynamics: uniqueness of the wave speeds
in the interior},
J.Differential Equations
162 (2000), 300-325.
%

\bibitem{[S]}{J.-E.Santos: }{\it Elastic wave propagation in fluid-saturated
porous media. I. The existence and uniqueness theorems}. RAIRO Mod\`el.
Math. Anal. Num\'er. 20 (1986), no. 1,
113-128.
%
\bibitem{[T]}{D.Tataru: }{\it Carleman estimates and
unique continuation for solutions to boundary value
problems},
J.Math.Pures Appl. 75 (1996), 367-408.

\bibitem{[Y]}{M.Yamamoto:}
{\it Carleman estimates for parabolic equations and
applications}, Inverse Problems 25 (2009) 123013
(75pp).



\end{thebibliography}
\end{document}